\documentclass{amsart}
\usepackage{amsmath,amsthm}
\usepackage {latexsym}
\usepackage{amssymb}
\usepackage[perpage,symbol]{footmisc}

\newcommand{\beal}{\begin{align}}
\newcommand{\enal}{\end{align}}
\newcommand{\bealn}{\begin{align*}}
\newcommand{\enaln}{\end{align*}}
\newcommand{\bear}{\begin{eqnarray}}
\newcommand{\eear}{\end{eqnarray}}
\newcommand{\beeq}{\begin{equation}}
\newcommand{\eneq}{\end{equation}}

\newcommand{\eps}{{\varepsilon}}
\newcommand{\R}{{\mathbb R}}

\newcommand{\calN}{{\mathcal N}}

\newcommand{\tileps}{{\tilde{\eps}}}

\def\bm{\left[ \begin{array}{cc}}
\def\endm{\end{array}\right]}

\def\eps{\varepsilon}

\def\bm{\left[\begin{matrix} }
\def\endm{\end{matrix}\right]}

\def\R{{\mathbb R}}

\newtheorem{theorem}{Theorem}
\newtheorem{lemma}[theorem]{Lemma}
\newtheorem{defi}[theorem]{Definition}

\newtheorem{proposition}[theorem]{Proposition}

\theoremstyle{remark}
\newtheorem{remark}[theorem]{Remark}

\renewcommand{\hat}{\widehat}
\renewcommand{\epsilon}{\eps}
\renewcommand{\tilde}{\widetilde}
\numberwithin{equation}{section}
\numberwithin{theorem}{section}

\begin{document}

\title[Optimal blow up co-rotaional WM]{Full blow-up range
  for co-rotaional wave maps to surfaces of revolution }

\author{Can Gao}



\thanks{}

\begin{abstract}
We construct blow-up solutions of the energy critical wave
map equation on $\R^{2+1}\to \mathcal N$ with polynomial blow-up rate
($t^{-1-\nu}$ for blow-up at $t=0$) in the case when $\calN$ is a surface of
revolution. Here we extend the blow-up range found by Carstea
($\nu>\frac 12$) based on the work by Krieger, Schlag and Tataru to
$\nu>0$. This work relies on and generalizes the recent result of
Krieger and the author where the target manifold is chosen as the standard
sphere.
\end{abstract}

\maketitle

\section{Introduction}
A \emph{wave map} is a map $u$ from $n+1$ dimensional
Minkowski space-time with signature $(-1,1,...,1)$ to a Riemannian
Manifold $\mathcal{N}$. It is defined as a critical point of the action functional, which is the following
Lagrangian
\[
\mathcal{L}(u): = \int_{\R^{2+1}}\langle\partial_{\alpha}u, \partial^{\alpha}u\rangle_{\mathcal N}\,d\sigma,\,\partial^{\alpha} = m^{\alpha\beta}\partial_{\beta}
\]
 where $\alpha = 0,1,...,n$, and
 $m^{\alpha\beta}$ is the Minkowski metric. 

The wave map $u: \mathbb{R}^{3+1}\to S^3$ has application to
  \emph{nonlinear sigma model}\cite{GL} from quantum field theory in modern
  physics, so it is very interesting to study the cases when target
  manifolds are spheres. The case $u: \mathbb{R}^{2+1} \to H^2$ is a model problem arising from the
  study of \emph{Einstein's equation}\cite{CM}. The curvature of the
  target manifold plays an important role in the global well-posedness
  properties of the corresponding equation. In the energy critical case (we will explain below what is
  energy critical) global well-posedness fails for the $S^2$ target,
  while it holds for $H^2$ (see below theorem \ref{thm:ST}
  and see\cite{K1,K2} and references therein). Another important observation
  is wave maps are the natural hyperbolic analogues of the much studied \emph{harmonic
map} heat flow, which in local coordinates is described by
$$\partial_tu^i=\Delta u^i+\sum_{\alpha=1}^n\Gamma^i_{jk}\partial_{\alpha}u^j\partial^{\alpha}u^k$$

Consider the following model equation
\begin{align}\label{model}
\Box u=N(u, \nabla u),\quad  (u, \partial_t
  u)|_{t=0}=(u_0, u_1)
\end{align}
for some smooth N(., .). Wave maps in local coordinates fall into this
category. Major studies of this problem fall into the following directions: i) local existence theory(strong local well-posedness);
ii) small data global existence theory(weak global well posed-ness);
iii) approaching the large data problem in the critical dimension
n=2 and hyperbolic target; iv) imposing symmetry: radial and equivariant wave maps in the case
n=2; v) singularity formation in the critical dimension. For details
of the main results in those directions, we refer the
reader to a very well-written survey paper on wave maps by Krieger \cite{K}
and the references therein.

In this paper, we study the blow-up solutions of energy critical co-rotational wave
map equation on $\R^{2+1}\to \mathcal N$ with polynomial blow-up rate in the case when $\calN$ is a surface of
revolution. Before we move further, we shall explain first about energy critical and definition of co-rotational.

\emph{Scaling constraints.} Assume that the set of solutions $u(t, x)$ of $(\ref{model})$ is invariant under the scaling
transformation $u(t, x)\to \lambda^{\alpha}u(\lambda t,\lambda x)$.Then
one introduces the \emph{critical Sobolev index} $s_c=\frac n2-\alpha$. Observe that the norm
$$\|u_0\|_{\dot{H}^{s_c}}+\|u_1\|_{\dot{H}^{s_c-1}}$$
is left invariant under the re-scaling. Note that $$s_c=\frac n2$$ for wave maps in the local coordinate formulation.

\emph{Energy constraints.} A quantity $$E[u] \gtrsim \|u\|_{H^{s_0}}+
\|u_t\|_{H^{s_0-1}}$$ which is preserved under
the flow. Then one distinguishes between: i) energy subcritical $s_c<s_0$: one expects global well-posedness, provided strong
local well-posedness in the full subcritical range, or also just for
some $s_c<s<s_0$; ii) energy critical $s_c=s_0$: global well-posedness hinges on
  fine structure of equation; iii) energy supercritical $s_c>s_0$: no global well-posedness for
  generic large data expected.

Note that when the background is $2+1$-dimensional, wave maps are {\it{energy
     critical}}. This means explicitly the following quantity 
\begin{equation}\label{eq:energy}
 \mathcal{E}(u): = \int_{\R^2}\big[|u_t|^2 + |\nabla_x u|^2|\big]\,dx
 \end{equation}
 is invariant under the intrinsic scaling (recall that $s_c=n/2$ in
 the local coordinate formulation)
$$u(t, x)\rightarrow u(\lambda t, \lambda x)
$$

\emph{Co-rotational wave maps.} A wave map $u: \mathbb{R}^{2+1}\to M$ is called \emph{equivariant} provided we have
$$u(t, \omega x) = \rho(\omega)u(t, x), \forall \omega\in S^1$$
Here $\rho(\omega)$
acts as an isometry on M and $\omega\in S^1$ acts on $\mathbb{R}^2$ in the canonical fashion as rotations. For
global well-posedness of equivariant wave maps we have the following important results by
Shatah, Tahvildar-Zadeh \cite{ST}
\begin{theorem}[Shatah, Tahvildar-Zadeh]\label{thm:ST}
Let the target $(M, g)$ be a warped product manifold satisfying a
suitable geodesic convexity condition. Then equivariant wave maps $u: \mathbb{R}^{2+1}\to M$
with smooth data stay globally regular.
\end{theorem}
However, the case $u:\mathbb{R}^{2+1}\to S^2$ does not satisfy the
hypotheses of the preceding theorem. Thus the discovery of the singularity for this
case is very crucial. We let $S^1$ act on $S^2$ by
means of rotations around the z-axis via $\rho(\omega)=k\omega, k\in
Z/ \{0\}$, $\omega\in S^1$. Fixing
a $k$, the wave map is then determined in terms of the polar angle, and becomes a
scalar equation on $\mathbb{R}^{1+1}$ as follows:
\begin{align}\label{sph-wave}
&-u_{tt}+u_{rr}+\frac 1{r}u_r=k^2\frac{\sin (2u)}{2r^2}
\end{align}
The case k = 1 in particular is called \emph{co-rotational}.

M. Struwe's fundamental work \cite{Struwe1} on the structure of singularities of
co-rotational Wave
maps shows that 
\begin{theorem}[Struwe]
Let $u$ be a smooth co-rotational wave map which cannot be
  smoothly extended past time $T$, there exists $t_i\rightarrow T,\quad \lambda_i\rightarrow +
  \infty$ s.t. on each fixed time slice $t = t_i$, we can write
$$
u(t_i, x) = Q(\lambda(t_i)x) + \epsilon(t_i, x)
$$
 where $Q$ is ground state (harmonic map) $Q: \mathbb{R}^2\rightarrow
 S^2$, while the local energy of $\epsilon$ converges to 0.
\end{theorem}
Furthermore, Struwe established an upper bound on the blow up rate
\begin{align}\label{blowuprate}
 \lim_{i\rightarrow\infty}\lambda(t_i)(T-t_i) = +\infty
\end{align}

The approach we take starts from \cite{KST0}, where the
 authors demonstrated a method of building finite time blow-up solutions for
 critical wave maps by adding corrections to an ansatz generated by
 rescaling the ground-state harmonic
 map to form an approximate solution and controlling the errors to
 zero. The blow-up rate from their paper is $\lambda(t)=t^{-1-\nu}$,
   with a blow-up range $\nu>\frac 12$. According to the work
   \cite{Struwe1} by M. Struwe (see above), this result is not
optimal (Replace $\lambda(t)=t^{-1-\nu}$ in $(\ref{blowuprate})$, one can see
that the optimal range for $\nu$ shall be $(0,\infty)$). 

In a joint work by the author and Krieger in \cite{CK}, the blow-up range is
extended to the full range $\nu>0$ which is optimal. 
It is also interesting to
consider the same problem in a more general situation when the target
manifold is a surface of revolution. A work on this case which is
parallel of \cite{KST0} was due to C\^arstea \cite{catalin}. However, as in
\cite{KST0}, the blow-up range in \cite{catalin} is not optimal. In
this paper, we will indicate how to combine the techniques of
\cite{catalin,CK} to obtain the optimal blow-up range in this
setting. For more detailed references concerning the blow-up dynamic of wave maps one can refer to \cite{CK}.

Let $\mathcal N$ be a surface of revolution equipped with a
Riemannian metric
$$ds^2=d\rho^2+g(\rho)^2d\theta$$
for $\mathcal N$ being produced by rotating the graph of a function $y=f(z)$ around the
$z$-axis. 
\begin{remark}
A detailed discussion of what properties $g$ shall
satisfy can be found in \cite{catalin}. Those properties will give the
relevant properties of the ground state (harmonic map) which one needs to use when proving some
intermediate conclusions when building the approximate
solutions. What this paper will focus on is the main difference and
changes raised because of the new setting of target manifold we
have. However, no changes are required according to the parts
of proofs relevant to $g$. Thus, we refer the
reader to \cite{catalin} for the details about what properties $g$
need to satisfy. 
\end{remark}

In the case of surfaces of revolution, the
equation for co-rotational wave maps takes a form similar to $(\ref{sph-wave})$. A simple computation (see \cite{catalin}) gives
\begin{align}\label{NW}
-\partial^2_tu+\partial_r^2u+\frac 1r\partial_ru=\frac{f(u)}{r^2},
\quad f(u)=g(u)g'(u).
\end{align}

Pick a stationary solution with finite energy for $(\ref{NW})$ as was shown in
\cite{catalin}. We state our result
\begin{theorem}\label{thm:Main} For any $\nu>0$, there exist $T>0$ and co-rotational initial data $(f, g)$ with 
\[
(f-\pi, g)\in H_{\R^2}^{1+\frac{\nu}{2}-}\times H_{\R^2}^{\frac{\nu}{2}-}
\]
 a\footnote{Here we use the identification of the wave map with a function $u(t, r)$ as before.} solution $u(t, r)$, $t\in (0, T]$ which blows up at time $t = 0$ and has the following representation: 
\[
u(t, r) = Q(\lambda(t)r) + \eps(t, r)
\]
where $\lambda(t) = t^{-1-\nu}$, and such that the function 
\[
(\theta, r) \longrightarrow \big(e^{i\theta}\eps(t, r), e^{i\theta}\eps_t(t, r)\big)\in H^{1+\nu-}(\R^2)\times H^{\nu-}(\R^2)
\]
uniformly in $t$. Also, we have the asymptotic as $t\rightarrow 0$
\[
\mathcal{E}_{loc}\big(\eps(t, \cdot)\big)\lesssim t^{\nu}\log^2 t
\]
\end{theorem} 

\section{A overview of the proof for theorem $\ref{thm:Main}$}
In the work on co-rotational wave maps to $S^2$ target by Krieger,
Schlag, and Tataru \cite{KST0}, it was found that solutions exist with the
blow-up rate $\lambda(t) = t^{-1-\nu}$, for the continnum of blow-up
rates of any $\nu > 1/2$. In a joint work of the author and Krieger
\cite{CK}, this range was extended to $\nu > 0$. Since the construction to
be described in this paper is based heavily on that of the previously
mentioned works, we recall for the convenience of the readers the
basic scheme. 

The method of construction relies on
building approximate solutions starting from the initial guess
$u(t,r) \approx Q(\lambda(t)r)$ where $Q(r)$ is the stationary ground
state. If one naively plugs in $Q(\lambda(t)r)$ into the equation, the
error term generated is $(r\lambda'(t))^2 Q''(\lambda(t)r) +
r\lambda''(t) Q'(\lambda(t)r)$, which turns out to be ``large''.
Thus one cannot directly use perturbative techniques to find the
solution. Instead, we first correct the error (within the past light
cone from the singularity) using an iterative scheme, until the error
becomes sufficiently small. In the following we will using the
notation $R = \lambda(t) r$."

\begin{theorem}\label{thm:approxsol} Assume $k\in \mathbf{N}$. There
  exists an approximate solution $u_{2k-1}(R)$  within the backwards light cone
from the singularity for \eqref{NW} which can be written as
\[
u_{2k-1}(t, r) = Q(R) + \frac{c_k}{(t\lambda)^2}R\log(1+R^2) +\frac{\tilde{c}_k}{(t\lambda)^2}R  + O\big(\frac{(\log(1+R^2))^2}{(t\lambda)^2}\big)
\]
with a corresponding error of size
\begin{align*}
e_{2k-1}
:&=\Big(-\partial_t^2+\partial^2_r+\frac{1}{r}\partial_r\Big)u_{2k-1}-\frac{f(u_{2k-1})}{2r^2}\\
&=(1-\frac{R}{\lambda t})^{-\frac{1}{2}+\nu}O\big(\frac{R(\log(1+R^2))^2}{(t\lambda)^{2k}}\big)
\end{align*}
Here the implied constant in the $O(\ldots)$ symbols are uniform in $t\in (0,\delta]$ for some $\delta = \delta(k)>0$ sufficiently small.
\end{theorem}

This is proved by means of an iterative scheme (see section
$\ref{sect:approx}$) that improves the error at each double
step. Actually at each step we approximately solve the wave equation
first close to $r=0$ then close to the light cone $r=t$. In both cases
it will reduce to solve an ODE (a Sturm-Louville equation). It is
important to observe here that the restriction $\nu>\frac 12$ imposed in \cite{catalin} does not come in at this stage; in fact, any $\nu>0$
will suffice. For the sake of readability,
only theorem $\ref{thm:approxsol}$ as well as the finer representation
of the errors as specified in $(\ref{eq:e_2k-1})$ will be used in the
final proof of the main theorem (the exact solution) in section
$\ref{sect:exact}$. The reader can treat section $\ref{sect:approx}$
as a black box if desired only up to these statements.

In section $\ref{sect:exact}$, we complete the approximate solution to
the exact one by adding correction via the ansatz
$u(t,r)=u_{2k-1}(t,r)+\epsilon(t,r)$. Before giving the relevant PDE
of such term $\epsilon$. We first renormalize the time $t$ into $\tau: = \nu^{-1} t^{-\nu}$, note that with respect to this time, we get
\[
\lambda(\tau) : = \lambda(t(\tau)) = (\nu\tau)^{\frac{1+\nu}{\nu}}
\]
We also have the re-scaled variable $R = \lambda(\tau)r$
respectively. We shall assume that
\[
|e_{2k-1}(t, r)|\lesssim \tau^{-N},\,r\leq t
\]
for some sufficiently large $N$, which is possible if we choose $k$
large enough. We shall also assume the fine structure of $e_{2k-1}$ as
in section $\ref{sect:approx}$, and more specifically as in \eqref{eq:e_2k-1}.
We can complete the approximate solution $u_{2k-1}$ to an exact
solution $u = u_{2k-1} + \eps$. , where $\eps$ solves the following
equation:
\begin{align}\label{eq:tilepsmain0}
&-\Big[\Big(\partial_{\tau}+\frac{\lambda_{\tau}}{\lambda}R\partial_R\Big)^2+\frac{\lambda_{\tau}}{\lambda}\Big(\partial_{\tau}+\frac{\lambda_{\tau}}{\lambda}R\partial_R\Big)\Big]\eps+\Big(\partial_R^2+\frac
1R\partial_R-\frac{f'(Q(R))}{R^2}\Big)\eps\nonumber\\
&=-\frac{1}{\lambda^2}[e_{2k-1}+N_{2k-1}(\eps)],
\end{align}
where
\begin{align}
N_{2k-1}(\eps)=\frac{1}{r^2}[f'(u_0)\eps-f(u_{2k-2}+\eps)+f(u_{2k-2})].
\end{align}
After changing of function $\tilde{\eps}(\tau,R)=R^{1/2}\eps(\tau,R)$,
$(\ref{eq:tilepsmain0})$ becomes
\begin{equation}\label{eq:tilepsmain1}
\big(-(\partial_{\tau} + \frac{\lambda_{\tau}}{\lambda}R\partial_R)^2 + \frac{1}{4}(\frac{\lambda_{\tau}}{\lambda})^2 + \frac{1}{2}\partial_{\tau}(\frac{\lambda_{\tau}}{\lambda})\big)\tilde{\eps} - \mathcal{L}\tilde{\eps} = \lambda^{-2}R^{\frac{1}{2}}\big(N_{2k-1}(R^{-\frac{1}{2}}\tilde{\eps}) + e_{2k-1}\big)
\end{equation}
The strategy is to formulate this equation in terms of the Fourier
coefficients of $\tilde{\eps}$ with respect to the generalized Fourier
basis associated with $ \mathcal{L}$ given by
\[
\mathcal{L} =  -\partial_R^2 + \frac{3}{4R^2}+V(R), \quad V(R)=-\frac{1}{R^2}[1-f'(Q(R))] 
\]
with $Q(R)$ the ground state. Dealing with $(\ref{eq:tilepsmain1})$,
one needs to develop some rather sophisticated spectral theory. The spectral theory of
$\mathcal{L}$ follows from \cite{catalin} (more exactly
\cite{KST0}), we refer the reader to \cite{KST0} to see a detailed
discussion. To find $\tilde{\eps}$, one employes a
fixed point argument in suitable Banach spaces, and it is here, in the
treatment of the nonlinear terms with singular weights, that the
restriction on $\nu$ comes in (see \cite{catalin,KST0}). More
precisely (see lemma 7.2 in \cite{catalin}), this condition is needed
there to make sufficient embedding between suitable function spaces to
control the nonlinear terms. 

In \cite{CK}, the authors overcome this
restriction (in the case while target manifold is sphere). We will
employ this method in our problem (while target manifold is surfaces of
revolution) in section \ref{sect:exact} which is as
following: 

Firstly, by a more closely analysis of the ’zeroth iterate’ (to be explained below) for $\tilde{\eps}$. We show that one can split this into the sum of two terms,
one of which has a regularity gain which lands us in the regime in \cite{KST0} is applicable, the other of which does not gain regularity
but satisfies an a priori $L^{\infty}$ bound near the symmetry axis $R=0$. So the relevant terms with a singular weight $R^{-3/2}$ at
$R=0$, such as $R^{-3/2}\tilde\epsilon^2$ (see section \ref{sect:exact})
can be estimated without adding any conditions for the regularity.
The reason why they can control the part of the zeroth iterate near
$R=0$ comes from the fact that the singular behavior of the approximate
solution from the first part of the construction and the error it generates is
localized to the boundary of the light cone. Then, by writing the equation for the distorted Fourier transform
of $\tilde{\eps}$ we will show that the higher iterates all differ from the zeroth iterate by terms with
a smoothness gain. This will then suffice to show the desired convergence.

\begin{remark}
The proof of Theorem $\ref{thm:Main}$, unsurprisingly, has large overlap with the
constructions of \cite{catalin,CK}. For brevity we will only indicate in
this note the modifications necessary, and will refer the reader to
\cite{catalin,CK} for the proofs of many intermediate steps.
\end{remark}
\begin{remark}
In the new situation, the main difficulty for proof of Theorem $\ref{thm:Main}$ is that we can not
write the nonlinear term explicitly. Thus in the relevant step (see
step 3 below) when
constructing the approximate solutions and in
the second part where the `perturbative scheme' is introduced for the
exact solutions, one needs to redo or adjust the proofs for the new
nonlinear source term. In \cite{CK}, the authors correct the
inaccuracies in \cite{KST0} according to the approximate solution step such as
the omission of some logarithm factors in the algebra of the special
function spaces. In out paper here, the different function spaces are
used correspondingly to fix such inaccuracies in \cite{catalin}. So some part of the arguments
need to be restated during the
construction of the approximate solutions.
\end{remark}

\section{Construction of the exact solutions}\label{sect:exact}
This is the very end of the proof of the main theorem. However this is
where the `key structure' is introduced following \cite{CK} to
make it possible to relax the constraint on $\nu$. For the readers who are interested in
the construction of the approximate solutions, we give the proof in
section \ref{sect:approx}. 

On the base that an approximate
solution has been constructed with a
corresponding error term which decays rapidly in the
renormalized time $\tau: = \nu^{-1} t^{-\nu}$, 
we can complete the approximate solution $u_{2k-1}$ to an exact
solution $u = u_{2k-1} + \eps$. After changing of function (which gives
us a new relevant $\tilde\eps$, see section 2) and applying a {\it distorted Fourier
transform}\footnote{Here the distorted Fourier transform is defined
  via combining one function $\phi(r,z)$ from the fundamental system
  for $\mathcal{L}-z$ and its inverse is given using the
  density function $\rho(\xi)$ of the spectral measure of $\mathcal{L}$, where
  $\mathcal{L}$ is a key operator raised from the exact solution's
  equation and $z\in\mathbb{C}$.\\ More precisely, the distorted
  Fourier transform is $$\mathcal{F}: \quad
  \hat{h}(\xi):=\int_0^{\infty}\phi(r,\xi)h(r)dr$$ when the inverse
  is $$\mathcal{F}^{-1}:\quad
  h(r):=\int_0^{\infty}\phi(r,\xi)\hat{h}(\xi)\rho(\xi)d\xi.$$ The
  detailed explanation for $\phi(r,z)$ and $\rho(\xi)$ is in
  \cite{CK,KST0}.} to the equation of $\tilde\eps$ (
($\ref{eq:tilepsmain1}$) in section 2):
\begin{equation}\label{eq:tilepsmain}
\big(-(\partial_{\tau} + \frac{\lambda_{\tau}}{\lambda}R\partial_R)^2 + \frac{1}{4}(\frac{\lambda_{\tau}}{\lambda})^2 + \frac{1}{2}\partial_{\tau}(\frac{\lambda_{\tau}}{\lambda})\big)\tilde{\eps} - \mathcal{L}\tilde{\eps} = \lambda^{-2}R^{\frac{1}{2}}\big(N_{2k-1}(R^{-\frac{1}{2}}\tilde{\eps}) + e_{2k-1}\big)
\end{equation}
One shall get a
equation of the Fourier coefficients, which we call the
\emph{transport equation}. 

The main difficulty is caused by the
operator $R\partial_R$ which is not diagonal in the Fourier basis. To
deal with this, we replace the distorted Fourier transform of $R\partial_Ru$
with $2\xi\partial_{\xi}$ modulo an error which will be treated perturbatively.  We define the error operator $\mathcal K$ by 
$$\widehat{R\partial_Ru}=-2\xi\partial_{\xi}\widehat{u}+\mathcal{K}\widehat{u}$$
where $\widehat{f}=\mathcal{F}f$ is the distorted Fourier
transform. 

To proceed further, we have to precisely understand the structure of the 'transference operator' $\mathcal{K}$. Make the
\begin{defi} We call an operator $\tilde{\mathcal{K}}$ to be 'smoothing', provided it enjoys the mapping property
\[
\tilde{\mathcal{K}}: L^{2, \alpha}_\rho\longrightarrow  L^{2, \alpha+\frac{1}{2}}_\rho\,\,\,\forall \alpha
\]
\end{defi}
For the definition of a weighted $L^2$-space $ L^{2,
  \alpha}_\rho$, we have
\[
\|u\|_{L^{2, \alpha}_\rho}: = \big(\int_0^\infty|u(\xi)|^2\langle\xi\rangle^{2\alpha}\rho(\xi)\,d\xi\big)^{\frac{1}{2}}
\]

If we put the terms with a
`smooth' property to the right hand
side of the equality in the transport equation. Then the Fourier coefficients (we call them
$x({\tau,\xi})$) of $\tilde\eps$ with respect to the generalized Fourier
basis satisfy
\begin{equation}\label{eq:Fourier4}
\mathcal{D}_{\tau}^2x + \xi x = f(x,\tilde{\eps}),
\end{equation}
where we have
the operator
\[
\mathcal{D}_{\tau}: = \partial_{\tau} - \frac{\lambda_{\tau}}{\lambda}[2\xi\partial_{\xi}+\frac{3}{2}+\frac{\rho'(\xi)\xi}{\rho(\xi)}]
\]
and
\begin{equation}\label{eq:Fourier5}\begin{split}
-f = &2\frac{\lambda_{\tau}}{\lambda}\mathcal{K}_0\big(\partial_{\tau} - \frac{\lambda_{\tau}}{\lambda}2\xi\partial_{\xi}\big)x + (\frac{\lambda_{\tau}}{\lambda})^2\big[\mathcal{K}^2 - (\mathcal{K} - \mathcal{K}_0)^2 - 2[\xi\partial_{\xi}, \mathcal{K}_0]\big]x\\
&+\partial_{\tau}(\frac{\lambda_{\tau}}{\lambda})\mathcal{K}_0x+\lambda^{-2}\mathcal{F}\big[R^{\frac{1}{2}}\big(N_{2k-1}(R^{-\frac{1}{2}}\tilde{\eps}) + e_{2k-1}\big)\big] - c\tau^{-2}x
\end{split}\end{equation}
For $\mathcal{K}_0$, according to \cite{catalin} we give it
as (see theorem 5.1\cite{catalin})
$$\mathcal{K}=-\Big(\frac 32 +\frac{\eta\rho'(\eta)}{\rho(\eta)}\Big)\delta_0(\xi-\eta)+\mathcal{K}_0.$$

\begin{remark}
Although the problem dealt in \cite{CK} is different than ours, the process at this stage is very close. We refer the readers to
\cite{CK} for those technical details we omit here when deducing the final transport
equation (mainly the straightforward computation) and below
for brevity.
\end{remark}

The explicit solution of \eqref{eq:Fourier4} is given as:  
\begin{lemma}[\cite{CK}] The equation \eqref{eq:Fourier4} is formally solved by the following parametrix
\begin{equation}\label{eq:Sparametrix}
x(\tau, \xi) = \int_{\tau}^\infty \frac{\lambda^{\frac{3}{2}}(\tau)}{\lambda^{\frac{3}{2}}(\sigma)}\frac{\rho^{\frac{1}{2}}(\frac{\lambda^2(\tau)}{\lambda^2(\sigma)}\xi)}{\rho^{\frac{1}{2}}(\xi)}S(\tau, \sigma, \lambda^2(\tau)\xi)f(\sigma, \frac{\lambda^2(\tau)}{\lambda^2(\sigma)}\xi)\,d\sigma =:(Uf)(\tau, \xi)
\end{equation}
\end{lemma}
One key fact from \cite{CK} is we have the following mapping property
of the parametrix with respect to suitable Banach spaces:
\begin{lemma}[lemma 5.6, \cite{CK}]\label{lem:ParaBounds}Introducing the norm
\[
\|f\|_{L_{\rho}^{2,\alpha;N}}: = \sup_{\tau>\tau_0}\tau^{N}\|f(\tau, \cdot)\|_{L^{2,\alpha}_{\rho}},
\]
we have
\[
\|Uf\|_{L_{\rho}^{2,\alpha+\frac{1}{2};N-2}}\lesssim \|f\|_{L_{\rho}^{2,\alpha;N}}
\]
provided $N$ is sufficiently large.
\end{lemma}

For the future reference, we will use the following norm:
\[
\|h\|_{H^{\alpha}_{\rho}}: = \big(\int_0^\infty x^2(\xi)\langle\xi\rangle^{2\alpha}\,\rho(\xi)\,d\xi\big)^{\frac{1}{2}}
\]
where $$h(R) = \int_0^\infty \phi(R, \xi)x(\xi)\rho(\xi)\,d\xi.$$

\subsection{Zeroth, first and higher iterative schemes}
After formulating \eqref{eq:Fourier4} 
 as an integral equation, we need to find a suitable fixed point, which will be
the desired $x(\tau, \xi)$. We construct these via
\begin{equation}\label{eq:fixedpoint}
x(\tau, \xi) = (Uf)(\tau, \xi)
\end{equation}
with $f(x,\tilde\eps)$ as in $(\ref{eq:Sparametrix})$. To find such a fixed point, we use the iterative scheme
\[
x_j(\tau, \xi) = (Uf_{j-1})(\tau, \xi),\,j\geq 1
\]
The function $f_j$ is
given as
\begin{align}\label{fj}
-f_j = &2\frac{\lambda_{\tau}}{\lambda}\mathcal{K}_0\big(\partial_{\tau} - \frac{\lambda_{\tau}}{\lambda}2\xi\partial_{\xi}\big)x_j + (\frac{\lambda_{\tau}}{\lambda})^2\big[\mathcal{K}^2 - (\mathcal{K} - \mathcal{K}_0)^2 - 2[\xi\partial_{\xi}, \mathcal{K}_0]\big]x_j\\
&+\partial_{\tau}(\frac{\lambda_{\tau}}{\lambda})\mathcal{K}_0x_j+\lambda^{-2}\mathcal{F}\big[R^{\frac{1}{2}}\big(N_{2k-1}(R^{-\frac{1}{2}}\tilde{\eps}_j) + \tilde{e}_{2k-1}\big)\big] - c\tau^{-2}x_j\nonumber
\end{align}

The zeroth iterate in turn is defined via
\[
x_0(\tau, \xi) = (U\lambda^{-2}\mathcal{F}\big[R^{\frac{1}{2}}\big(e_{2k-1}\big)\big])(\tau, \xi);
\]
We have the following proposition proved in \cite{CK}
\begin{proposition}[proposition 5.7, \cite{CK}]
Replacing $e_{2k-1}$ with $\tilde{e}_{2k-1}\in H^{\frac
  {\nu}{2}-}_{RdR}$ where $\tilde{e}_{2k-1}|_{r\leq t}=e_{2k-1}$, we can write $$x_0=x_0^{(1)}+x_0^{(2)}$$ where 
$$x^{(1)}\in\tau^{-N}L_{\rho}^{2,\frac 12+\frac{\nu}{2}-}, \quad
x^{(2)}\in\tau^{-N}L_{\rho}^{2,1+\frac{\nu}{2}-}$$
and
also $$\chi_{R<1}\tilde\epsilon^{(1)}_0(\tau,R)=\chi_{R<1}\int_{0}^{\infty}\phi(R,\xi)x^{(1)}_0(\tau,\xi)\rho(\xi)d\xi\in\tau^{-N}R^{\frac
32}L^{\infty},\quad \chi_{R\geq 1}|\tilde{\epsilon}_0^{(1)}|\lesssim \tau^{-N}$$
\end{proposition}

We can rephrase it as following, which is identical to Corollary
5.9 in \cite{CK}. 
\begin{proposition}\label{cor:refine}Denote by $P_{\lambda}$ the frequency localizers
\[
\mathcal{F}\big(P_{<\lambda} f\big)(\xi) = \chi_{<\lambda}(\xi)\big(\mathcal{F}f\big)(\xi)
\]
where $ \chi_{<\lambda}(\xi)$ is a smooth cutoff function localizing to $\xi\lesssim \lambda$, as in \cite{KST0}; here $\lambda$ is a dyadic number. Then we have
\[
\chi_{R<1}P_{<\lambda}\tileps_0^{(1)}\in \tau^{-N}R^{\frac{3}{2}}L^\infty
\]
uniformly in $\lambda>1$. Furthermore, for any integer $l\geq 0$, we have
\[
\nabla_R^l R^{-\frac{3}{2}}P_{<\lambda}\tileps_0^{(1)} = O(\tau^{-N})
\]
uniformly in $\lambda>1$.
\end{proposition}
\begin{remark}
This is the key structure from \cite{CK}, with which the we are able to invoke
lemma \ref{lm:NLS} to control the nonlinear term and prove
$(\ref{NLT})$ (see below).
\end{remark}

Based on lemma $\ref{lem:ParaBounds}$, we know $$\|Uf_{j-1}\|_{L_{\rho}^{2,1+\frac{\nu}{2}-}}\lesssim\|f_{j-1}\|_{L_{\rho}^{2,\frac 12+\frac{\nu}{2}-}}$$

For the first iterate, the estimate for the most terms in
$(\ref{fj})$ follows the same arguments in
\cite{CK}. We list the unchanged results
(see \cite{CK} for proof) as following
\begin{align*}
&(\partial_{\tau}-\frac{\lambda_{\tau}}{\lambda}2\xi\partial_{\xi})x_0\in
\tau^{-N-1}L_{\rho}^{2,\frac{\nu}{2}-}\\
&2\frac{\lambda_{\tau}}{\lambda}\mathcal{K}_0(\partial_{\tau}-\frac{\lambda_{\tau}}{\lambda}2\xi\partial_{\xi})x_0\in\tau^{-N-2}L_{\rho}^{2,\frac
  12+\frac{\nu}{2}-}\\
&(\frac{\lambda_{\tau}}{\lambda})^2[\mathcal{K}^2-(\mathcal{K}-\mathcal{K}_0)^2-2[\xi\partial_{\xi},\mathcal{K}_0]]x_0\in
\tau^{N-2}L_{\rho}^{2,1+\frac{\nu}{2}-}\\
&\partial_{\tau}(\frac{\lambda_{\tau}}{\lambda})\mathcal{K}_0x_0-c\tau^{-2}x_0\in
\tau^{-N-2}L_{\rho}^{2,\frac 12+\frac{\nu}{2}-}
\end{align*}
For the nonlinear term, which is the key of the whole argument, we will prove the following in the next
section (according to Lemma 3.4)
\begin{align}\label{NLT}
\lambda^{-2}R^{\frac 12}N_{2k-1}(R^{-\frac 12}\tilde
\epsilon)\in \tau^{-N-2}L_{\rho}^{2,\frac 12+\frac{\nu}{2}-}
\end{align}

Let us for now accept the facts above and conclude here the key conclusion in this step
\[
\big\|x_1(\tau, \cdot)-x_{0}(\tau, \cdot)\big\|_{L^{2,1+\frac{\nu}{2}-}_{\rho}}\lesssim N^{-1}\tau^{-N},
\]
\[
\big\|\big(\partial_{\tau} - \frac{\lambda_{\tau}}{\lambda}2\xi\partial_{\xi}\big)\big(x_1(\tau, \cdot)-x_{0}(\tau, \cdot)\big)\big\|_{L^{2,\frac{1}{2}+\frac{\nu}{2}-}_{\rho}}\lesssim N^{-1}\tau^{-N-1}
\]

Then we define $$\tilde\epsilon_1=\int_0^{\infty}\phi(R,\xi)\big(x_1(\tau, \cdot)-x_{0}(\tau, \cdot)\big)\rho(\xi)d\xi+\int_0^{\infty}\phi(R,\xi)
x_0(\tau,\xi)\rho(\xi)d\xi$$ which will allow us to write
$$\tilde\epsilon_1=\tilde\epsilon^{(1)}(\tau,\cdot)+\tilde\epsilon^{(2)}(\tau,\cdot)$$
$\tilde\epsilon^{(1)}(\tau,\cdot)$ and
$\tilde\epsilon^{(2)}(\tau,\cdot)$ satisfy exactly the kind of structure we need to invoke the bound
for nonlinear source term in lemma \ref{lm:NLS}. Continuing running the
iterate scheme will give us the bounds 
\[
\big\|x_j(\tau, \cdot)-x_{j-1}(\tau, \cdot)\big\|_{L^{2,1+\frac{\nu}{2}-}_{\rho}}\lesssim N^{-j}\tau^{-N},
\]
\[
\big\|\big(\partial_{\tau} - \frac{\lambda_{\tau}}{\lambda}2\xi\partial_{\xi}\big)\big(x_j(\tau, \cdot)-x_{j-1}(\tau, \cdot)\big)\big\|_{L^{2,\frac{1}{2}+\frac{\nu}{2}-}_{\rho}}\lesssim N^{-j}\tau^{-N-1}
\]
This will close the fix point argument which proves we have
$$x_{\tau,\xi}\in H^{\frac 12+\frac{\nu}{2}-},\quad \partial_{\tau}x_{\tau,\xi}\in H^{\frac{\nu}{2}-}.$$ Through lemma 7.1 in
\cite{catalin} (it was proven in \cite{KST0}):
\begin{lemma}\label{3.2}
Assume $|\alpha|<\frac{\nu}{2}+\frac 34$, $g\in IS(1,\mathcal{Q})$. Then we have
$$\|gf\|_{H^{\alpha}_{\rho}}\lesssim \|f\|_{H^{\alpha}_{\rho}}$$
\end{lemma}

It indicates the existence of the exact solution
$\eps(\tau,\cdot)\in \tau^{-N}H^{1+{\nu}-}_{\mathbb{R}^2}$,
as well as $\partial_{\tau}\epsilon(\tau,\cdot)\in
\tau^{-N-1}H^{{\nu}-}_{\mathbb{R}^2}$.

\subsection{The nonlinear source terms}
We will give an analysis
to the new nonlinear source term to complete our work in this section. We recall the following formula for the main source term:
\begin{align}
\lambda^{-2}R^{\frac 12}N_{2k-1}(R^{-\frac 12}\tilde
\epsilon)&=\frac{1}{R^2}\big[f'(u_0)\tilde\epsilon-f(u_{2k-2}+R^{-\frac
12}\tilde\epsilon)R^{\frac 12}+f(u_{2k-2})R^{\frac 12}\big]\\
&=\frac
{1}{R^2}\big[f'(u_0)-f'(u_{2k-2})\big]\tilde\epsilon-\frac{1}{R^{\frac
  32}}\sum_{l\geq 2}\frac{1}{l!}f^{(l)}(u_{2k-2})\big(R^{-\frac 12}\tilde\epsilon\big)^l
\end{align}
According to the preceding proposition, we have 
 \[
 x_0\in \tau^{-N}L^{2, \frac{1}{2}+\frac{\nu}{2}-}_{\rho}
 \]
 whence 
 \[
 \tileps_0(\tau, \cdot)\in \tau^{-N}H^{\frac{1}{2}+\frac{\nu}{2}-}_{\rho}
 \]
 This means that for the source terms, we need at least
 $H^{\frac{\nu}{2}-}_{\rho}$-regularity. In fact, we can do much
 better for the term $\frac
{1}{R^2}\big[f'(u_0)\tilde\epsilon-f'(u_{2k-2})\tilde\epsilon\big]$. Recall
that $$u_{2k-2}=u_0+\sum_{j=1}^{2k-2}v_j$$
where we have
$$v_{2k-1}\in \frac{1}{(t\lambda)^{2k}}IS^3\big(R(\log R)^{2k-1},
\mathcal{Q}_{k-1}\big), \quad v_{2k}\in \frac{1}{(t\lambda)^{2k+2}}IS^3\big(R^3(\log R)^{2k-1}, \mathcal{Q}_{k}\big)$$
which implies
$$u_{2k-2}-u_0\in \frac {1}{(t\lambda)^2}IS^3(R\log R,\mathcal Q)$$
Moreover, we recall some useful results in \cite{catalin,KST0}.
\begin{lemma}[lemma 3.9-10, \cite{catalin}]\label{catalinlm}
$f^{(2k)}(u_0)\in IS^1(R^{-1})$ and $f^{(2k+1)}(u_0)\in IS^0(1)$. Moreover, if $$z\in\frac{1}{(t\lambda)^2}IS^1(R(\log R), \mathcal Q),$$
then $$f^{(2k)}(u_0+z(R))\in \frac{1}{(t\lambda)^2}IS^1(R(\log R), \mathcal Q)$$
and $$f^{(2k+1)}(u_0+z(R))\in IS^0(1, \mathcal Q).$$
\end{lemma}
Thus for $\frac
{1}{R^2}\big[f'(u_0)\tilde\epsilon-f'(u_{2k-2})\tilde\epsilon\big]$, we have
$$\frac{f'(u_0)-f'(u_{2k-2})}{R^2}=\frac{\sum_{l\geq 2}\frac{1}{l!}f^{(l)}(u_0)(u_{2k-2}-u_0)^l}{R^2}\in \frac{1}{(t\lambda)^2}IS(1,\mathcal{Q})$$
and lemma \ref{3.2} will give us the following bound
\begin{align}\label{NLS:1}
\|\frac
{1}{R^2}\big[f'(u_0)\tilde\epsilon-f'(u_{2k-2})\tilde\epsilon\big]\|_{H^{\frac 12+\frac {\nu}{2}-}}\lesssim (t\lambda)^{-2}\|\tilde\epsilon\|_{H^{\frac 12+\frac {\nu}{2}-}}
\end{align}
To deal with the rest `truly' nonlinear terms, we first split them into two parts
\begin{align}
&\frac{1}{R^{\frac
  32}}\sum_{l\geq 2}\frac{1}{l!}f^{(l)}(u_{2k-2})\big(R^{-\frac 12}\tilde\epsilon\big)^l=\nonumber\\
&\frac{1}{R^{\frac
  32}}\sum_{l\geq 1}\frac{1}{l!}f^{(2l)}(u_{2k-2})\big(R^{-\frac 12}\tilde\epsilon\big)^{2l}\label{NLS:even}\\
&+\frac{1}{R^{\frac
  32}}\sum_{l\geq 1}\frac{1}{l!}f^{(2l+1)}(u_{2k-2})\big(R^{-\frac 12}\tilde\epsilon\big)^{2l+1}\label{NLS:odd}
\end{align}

We can write $(\ref{NLS:even})$ in the form
$$R^{-\frac 32}\tilde\epsilon^2\sum_{l\geq 1}\frac{1}{l!}\frac{f^{(2l)}(u_0+u_{2k-2}-u_0)}{R}\big(R^{-1}\tilde\epsilon^2\big)^{l-1}$$
and meanwhile write $(\ref{NLS:odd})$ as
$$R^{-3}\tilde\epsilon^3\sum_{l\geq 1}\frac{1}{l!}f^{(2l+1)}(u_0+u_{2k-2}-u_0)\big(R^{-1}\tilde\epsilon^2\big)^{l-1}$$
According to Lemma $\ref{catalinlm}$, we observe that
$$\frac{f^{(2l)}(u_0+u_{2k-2}-u_0)}{R},\quad f^{(2l+1)}(u_0+u_{2k-2}-u_0)
\in IS^0(1,Q).$$ Thus via Lemma
$(\ref{3.2})$, we can estimate the $H_{\rho}^{\alpha}$ norm of $(\ref{NLS:even})$
and $(\ref{NLS:odd})$ by the $H_{\rho}^{\alpha}$ norm of
$$R^{-\frac 32}\tilde{\epsilon}^2q(R^{-1}\tilde\epsilon^2),\quad R^{-3}\tilde\epsilon^3q(R^{-1}\tilde\epsilon^2)$$
where $\alpha$ here is $\frac 12+\frac {\nu}{2}-$ and $q(\cdot)$ is a
real analytic function.

We recall a very
technical and crucial
lemma proved in \cite{CK} 
\begin{lemma}[lemma 5.12, \cite{CK}]\label{lm:5.12}
Assume that all of $f, g, h$ are either in $H^{\frac{1}{2}+\frac{\nu}{2}-}_{\rho}\cap R^{\frac{3}{2}}L^\infty$ as well as with their frequency localized constituents $P_{<\lambda}(\cdot)\in \log \lambda R^{\frac{3}{2}}L^\infty$ and $\chi_{R<1}\nabla_R^l\big(R^{-\frac{3}{2}}P_{<\lambda}(\cdot)\big)\in L^\infty$, $l\geq 0$, uniformly in $\lambda>1$,  or in $H^{1+\frac{\nu}{2}-}_{\rho}$. Then we have
\[
R^{-3}fgh\in H^{\frac{1}{2}+\frac{\nu}{2}-}\cap R^{\frac{3}{2}}L^\infty,\,P_{<\lambda}(R^{-3}fgh)\in \log\lambda R^{\frac{3}{2}}L^\infty,\,P_{<\lambda}(R^{-3}fgh)\in RL^\infty
\]
with the latter two inclusions uniformly in $\lambda>1$. Also, if $h_j\in H^{\frac{1}{2}+\frac{\nu}{2}-}_{\rho}\cap R^{\frac{3}{2}}L^\infty$ and further $P_{<\lambda}h_j\in RL^\infty$ as well as $\chi_{R<1}\nabla_R^l\big(R^{-1}P_{<\lambda}h_j\big)\in L^\infty$, $l\geq 0$, uniformly in $\lambda$, or else $h_j\in H_{\rho}^{1+\frac{\nu}{2}-}$, for $j = 1,2,\ldots, 2N$, then we have
\[
R^{-3}fgh \prod_{j=1}^N(\frac{1}{R}h_{2j} h_{2j-1})\in H^{\frac{1}{2}+\frac{\nu}{2}-}
\]
We also get
\[
R^{-\frac{3}{2}}fg\prod_{j=1}^N(\frac{1}{R}h_{2j} h_{2j-1})\in H^{\frac{1}{2}+\frac{\nu}{2}-}
\]
\end{lemma}

Invoke the conclusion from lemma \ref{lm:5.12}, one can prove: 
\begin{lemma}\label{lm:NLS}
Providing
$$\|\tilde{\epsilon}\|_{H^{\frac 12+\frac{\nu}{2}-}\cap R^{\frac
    32}L^{\infty}}\lesssim 1,\quad \|R^{-\frac
  32}P_{<\lambda}\tilde{\epsilon}\|_{L^{\infty}}\lesssim 1,\quad
\|\chi_{R<1}\nabla^l_R(R^{-\frac 32}
P_{<\lambda}\tilde\epsilon)\|_{L^{\infty}}\lesssim 1$$
uniformly in $\lambda>1$ $l\geq 0$, we have
\begin{align*}
&\|\frac
{1}{R^2}\big[f'(u_0)\tilde\epsilon-f'(u_{2k-2})\tilde\epsilon\big]\|_{H_{\rho}^{\frac 12+\frac {\nu}{2}-}}\lesssim (t\lambda)^{-2}\|\tilde\epsilon\|_{H_{\rho}^{\frac 12+\frac {\nu}{2}-}}\\
&\|\frac{1}{R^{\frac
  32}}\sum_{l\geq 1}\frac{1}{l!}f^{(2l)}(u_{2k-2})\big(R^{-\frac 12}\tilde\epsilon\big)^{2l}\|_{H_{\rho}^{\frac 12+\frac {\nu}{2}-}}\lesssim \|\tilde\epsilon\|^2_{H_{\rho}^{\frac 12+\frac {\nu}{2}-}\cap R^{\frac 32}L^{\infty}}\\
&\|\frac{1}{R^{\frac
  32}}\sum_{l\geq 1}\frac{1}{l!}f^{(2l+1)}(u_{2k-2})\big(R^{-\frac 12}\tilde\epsilon\big)^{2l+1}\|_{H_{\rho}^{\frac 12+\frac {\nu}{2}-}}\lesssim \|\tilde\epsilon\|^3_{H_{\rho}^{\frac 12+\frac {\nu}{2}-}\cap R^{\frac 32}L^{\infty}}
\end{align*}
The last two estimates' right hand side space can be replaced by
$H^{\frac 12+\frac {\nu}{2}-}$ with a change of the bound of
$\tilde\epsilon$ by $$\|\tilde\epsilon\|_{H_{\rho}^{\frac
    12+\frac{\nu}{2}-}}\lesssim 1.$$
\end{lemma}

\section{The construction of the approximate solutions}\label{sect:approx} 

To build the approximate solution as in theorem $\ref{thm:approxsol}$, we follow the
scheme in \cite{KST0}. We start from the stationary harmonic
map\footnote{ The properties of ground state are needed to prove the spectral theory of $\mathcal{L}$. Since we will employ the same
  spectral theory as it is in \cite{catalin}, we refer the reader to section 2 \cite{catalin} for the
  discussion of properties of such ground states,} $Q(R)$. Setting $R =
\lambda(t) r$ we take $u_0(t,x) = Q(\lambda(t)x)$ for $\lambda(t) =
t^{-1-\nu}$ and then add corrections $v_k$
iteratively $u_k=u_0+\sum_{j=1}^kv_k$. In a first approximation we linearize the equation for the correction
$\epsilon=u-u_{k}$ around $\epsilon=0$ and substitute $u_{k}$ by
$u_0$. Then we have the linear approximate equation
$$\big(-\partial^2_t+\partial^2_r+\frac
1r \partial_r\big)\epsilon-\frac{1}{r^2}f'(u_0)\epsilon\thickapprox -e_k$$
From here we split into two different cases: considering the case $r\ll t$
when we expect the time derivative to play a lesser role thus we
neglect it (where $(\ref{equ:odd})$ below comes from); considering the case
$r\thickapprox t$ when the time and spatial derivative have the same
strength. We can identify another principal variable, namely $a=r/t$
and think of $\epsilon$ as a function of $\epsilon(t,a)$ so we can
reduce this case to a Strum-Liouville problem in $a$ which becomes
singular at $a=1$ (where$(\ref{equ:even})$ comes from). After each step of adding the
correction, we also estimate the size of the errors. This makes each
round of the scheme with four steps to go. For odd and even
steps, we have different equations for the corrections $v_k$:
\begin{align}
&\big(\partial^2_r+\frac
1r\partial_r-\frac{1}{r^2}f'(u_0)\big)v_{2k+1}=-e^0_{2k}\label{equ:odd}\\
&\big(-\partial_t^2+\partial^2_r+\frac 1r\partial_r-\frac {1}{r^2}\big)v_{2k+2}=-e^0_{2k+1}\label{equ:even}
\end{align}with Cauchy zero data\footnote{The coefficients are
  singular at $r=0$, therefore this has to be given a suitable
  interpretation below (see remark \ref{explain:singular}).} at $r=0$, and\footnote{There is a typo in
\cite{catalin} for the sign of the term $f(u_{2k-2})$. This
does not influence the result in \cite{catalin} but it matters for
our analysis for the nonlinear source terms in later section.}
where
\begin{align}
&e_k=\big(-\partial_t^2+\partial^2_r+\frac 1r\partial_r\big)u_k-\frac{1}{r^2}f(u_k)\\
&e_{2k+1}=e^1_{2k}-\partial^2_tv_{2k+1}+N_{2k+1}(v_{2k+1}),\quad
e_{2k}=e^1_{2k-1}+N_{2k}(v_{2k})\\
&N_{2k-1}(v)=\frac{1}{r^2}[f'(u_0)v-f(u_{2k-2}+v)+f(u_{2k-2})]\\
&N_{2k}(v)=\frac v{r^2}-\frac 1{r^2}[f(u_{2k-1}+v)-f(u_{2k-1})]
\end{align}
\begin{remark}
Note here a technical detail is we split $e_k$ into $e_k=e^0_k+e^1_k$ where $e^0_k$ is the
so-called principle part and the rest $e^1_k$, the so-called higher
order part, will be left and merge into the next step while analyzing
the error $v_{k+1}$ (will be precise below in step 1 and 3). Also we
will switch to the principle variable `a' for equation (\ref{equ:even}) in
step 3 as already mentioned in the above section.
\end{remark}
To formalize this scheme we need to define suitable function spaces in
the light-cone 
$$\mathcal{C}_0=\{(t,r):0\leq r< t, 0<t<t_0\}$$
to put our successive corrections and errors. They are
following closely from those in \cite{CK}.\footnote{One shall note that those definitions are
very natural according to a direct computation for the first round of
the iterative scheme (see \cite{CK} for the case when target manifold
is sphere). }
\begin{defi}\label{defi:Q_n}For $i\in \mathbf{N}$, let $j(i) = i$ if $\nu$ is irrational, respectively $j(i) = 2i^2$ if $\nu$ is rational.
Then
\begin{itemize}
\item $\mathcal{Q}$ is the algebra of continuous functions $q: [0,1]\rightarrow \R$ with the following properties:
\\
(i) $q$ is analytic in $[0,1)$ with even expansion around $a = 0$.\\
(ii) near $a = 1$ we have an absolutely convergent expansion of the form
\begin{align*}
q(a) = &q_0(a) + \sum_{i=1}^\infty (1-a)^{\beta(i)+\frac{1}{2}}\sum_{j=0}^{j(i)}q_{i,j}(a)\big(\log(1-a)\big)^j\\
&+ \sum_{i=1}^\infty (1-a)^{\tilde{\beta}(i)+\frac{1}{2}}\sum_{j=0}^{j(i)}\tilde{q}_{i,j}(a)\big(\log(1-a)\big)^j
\end{align*}
with analytic coefficients $q_0, q_{i,j}$, and $\beta(i) = i\nu$, $\tilde{\beta}(i) = \nu i+\frac{1}{2}$.
\item $\mathcal{Q}_n$ is the algebra which is defined similarly, but also requiring $q_{i,j}(1) = 0$ if $i\geq 2n+1$.
\end{itemize}
\end{defi}
We also define the space of functions obtained by differentiating $\mathcal{Q}_n$:
\begin{defi}\label{defi:Q'_n}
Define $\mathcal{Q}'$ as in the preceding definition but replacing $\beta(i)$ by $\beta'(i): = \beta(i) - 1$, and similarly for $\mathcal{Q}'_n$.
\end{defi}

\begin{defi}\label{defi:S^n} $S^n(R^k(\log R)^l)$ is the class of analytic functions $v: [0,\infty)\rightarrow \R$ with the following properties:\\
(i) $v$ vanishes of order $n$ at $R = 0$.\\
(ii) $v$ has a convergent expansion near $R = \infty$
\[
v = \sum_{\substack{0\leq j\leq l+i\\ i\geq 0}}c_{ij}R^{k-i}(\log R)^j
\]
\end{defi}
The final function
space $S^m(R^k(\log R)^l,Q_n)$ is defined slightly different than Definition 3.5 in \cite{CK}
where we add an extra $`b'$ into it. This is simply for applying the
results from \cite{catalin} later. We state it here precisely.
\begin{defi}\label{defi:S^mQ_n}(Definition 3.5, \cite{CK})  
Introduce the symbols
\[
b = \frac{\big(\log(1+R^2)\big)^2}{(t\lambda)^2},\,b_1 = \frac{\big(\log(1+R^2)\big)}{(t\lambda)^2},\,b_2 = \frac{1}{(t\lambda)^2}
\]
Pick $t$
  sufficiently small such that all $b,b_1, b_2$, when restricted to the light cone $r\leq t$ are of size at most $b_0$.
\begin{itemize}
\item $S^m(R^k(\log R)^l, \mathcal{Q}_n)$ is the class of analytic functions $v: [0,\infty)\times [0,1)\times [0,b_0]^3\rightarrow\R$ so that\\
(i) $v$ is analytic as a function of $R, b, b_1, b_2$,
\[
v: [0,\infty)\times [0, b_0]^3\rightarrow {\mathcal{Q}}_n
\]
(ii) $v$ vanishes to order $m$ at $R = 0$.\\
(iii) $v$ admits a convergent expansion at $R = \infty$,
\[
v(R,\cdot,b,b_1,b_2) = \sum_{\substack{0\leq j\leq l+i\\ i\geq 0}}c_{ij}(\cdot,b,b_1,b_2)R^{k-i}(\log R)^j
\]
where the coefficients $c_{ij}: [0, b_0]^3\rightarrow \mathcal{Q}_n$ are analytic with respect to $b,b_{1,2}$.
\item $IS^m(R^k(\log R)^l, \mathcal{Q}_n)$ is the class of analytic functions $w$ inside the cone $r<t$ which can be represented as
\[
w(t, r) = v(R, a,b, b_1, b_2),\,v\in S^m(R^k(\log R)^l, \mathcal{Q}_n)
\]
and $t>0$ sufficiently small.
\end{itemize}
\end{defi}
\begin{remark}\label{explain:singular}
The functional spaces $S^m(R^k(\log R)^l,Q_n)$ satisfy some good
asymptotic behaviors (for example, they vanish in order m at $R=0$) so the existence of the solutions to equation
$(\ref{equ:odd})$ and $(\ref{equ:even})$ will make sense in those spaces although
the coefficients are singular at $R=0$ in general.
\end{remark}

Following the method in \cite{KST0}, the idea for proving theorem $\ref{thm:approxsol}$ is to inductively
show that we can choose the corrections $v_k$ to be in relevant
function spaces:
\begin{equation}\label{eq:v_2k-1}
v_{2k-1}\in \frac{1}{(t\lambda)^{2k}}IS^3\big(R(\log R)^{2k-1}, \mathcal{Q}_{k-1}\big)
\end{equation}
\begin{equation}\label{eq:e_2k-1}
t^2e_{2k-1}\in \frac{1}{(t\lambda)^{2k}}IS^1\big(R(\log R)^{2k-1}, \mathcal{Q}'_{k-1}\big)
\end{equation}
\begin{equation}\label{eq:v_2k}
v_{2k}\in \frac{1}{(t\lambda)^{2k+2}}IS^3\big(R^3(\log R)^{2k-1}, \mathcal{Q}_{k}\big)
\end{equation}
\begin{equation}\label{eq:e_2k}
t^2e_{2k}\in \frac{1}{(t\lambda)^{2k}}\big[IS^1\big(R^{-1}(\log
R)^{2k}, \mathcal{Q}_{k}\big) + \langle b, b_1, b_2\rangle[IS^1\big(R(\log R)^{2k-1}, \mathcal{Q}'_{k}\big)\big]
\end{equation}
and the starting error $e_0$ satisfying
\[
e_0\in IS^1(R^{-1})
\]

Here we denote by $\langle b, b_1, b_2\rangle$ the ideal generated by $b,b_1,b_2$ inside the algebra generated
by $b,b_1,b_2$. Now we give a brief outline of the proof for $\ref{thm:Main}$:

\begin{proof}

First one shall check $e_0\in
IS^1(R^{-1})$, this can be done by a direct computation (see step 0 in
\cite{catalin}). Then assuming $(\ref{eq:v_2k-1}-\ref{eq:e_2k})$ hold
up to $k-1$, the first task would be proving
$(\ref{eq:v_2k-1})$ for $k$. \\

{\bf{Step 1}}: {\it{For $e_{2k-2},
    k\geq 1$, proves $v_{2k-1}$ satisfies \eqref{eq:v_2k-1}.}} \\

For this one first needs to choose the right
`principal part' of $e_{2k-2}$ which we call $e^0_{2k-2}$. This is
done by throwing away the
`higher order parts', which we call $e^1_{2k-2}$ and which belong to the same space
as $e_{2k-1}$. The way to do it is as following: when $k=1$ we let
$e^0_0:=e_0$, if $k>1$, we let $e^0_{2k-2}:=e_{2k-2}(R,a,0)$ with the
setting $b,b_1,b_2=0$. By changing
into variable $R$, equation $(\ref{equ:odd})$ becomes:
$$(t\lambda)^2Lv_{2k-1}=-t^2e^0_{2k-2}.$$
Here the operator $L$ is 
$$ L:= \partial^2_R+\frac{1}{R}\partial_R-\frac{f'(u_0)}{R^2}$$
 To get the
desired result, one needs to prove the following lemma:
\begin{lemma}
The solution of $Lv=\varphi\in S^1(R^{-1}(\log R)^{2k-2})$, with
$v(0)=v'(0)=0$, has the regularity
$$v\in S^3(R(\log R)^{2k-1}).$$
\end{lemma}

This is already proven as
Lemma 3.11 in
\cite{catalin}, so we conclude
$(\ref{eq:v_2k-1})$.\\

{\bf{Step 2}}: {\it{Choose $v_{2k-1}$ as in \eqref{eq:v_2k-1} with error $e_{2k-1}$ satisfying \eqref{eq:e_2k-1}.}}\\
 
According to the definition of $e_{2k-1}$
above, we have 
$$t^2e_{2k-1}=t^2e^1_{2k-2}-t^2\partial^2_tv_{2k-1}+t^2N_{2k-1}(v_{2k-1})$$
Since in the former step we treat $a$ as a parameter and now we
will defreeze it, some extra terms will show up while calculating the error
$e_{2k-1}$. To be more precise, the amended term
$t^2e_{2k-1}$ we need to deal with is as following (note that
$t^2e^1_{2k-2}$ is proved automatically thanks to the assumptions)
$$t^2e_{2k-1}=t^2N_{2k-1}(v_{2k-1})+E^tv_{2k-1}+E^av_{2k-1}$$
where $E^2v_{2k-1}$ is the term in $\partial^2_t v_{2k-1}$ with no
derivation on the $a$ variable, and the term $E^av_{2k-1}$ is the
terms in $(-\partial^2_t+\partial^2_r+\frac 1r\partial_r)v_{2k-1}$
where derivative hits the $a$ variable (the extra terms from defreezing of $a$ are
included here). To prove all those terms in $\ref{eq:e_2k-1}$, we refer the reader to step 2 in \cite{catalin}.\\

{\bf{Step 3}}: {\it{Given $e_{2k-1}$ as in \eqref{eq:e_2k-1}, construct $v_{2k}$ as in \eqref{eq:v_2k}}}\\

Here we have to diverge slightly from \cite{catalin}, since our
definition of the algebra $S^m(R^k\log R^l)$ is different (we follow
the definition in 
\cite{CK}). Since the equation $(\ref{equ:even})$ for $v_{2k}$ is
identical with equation $(3.2)$ for $v_{2k}$ in \cite{CK}. We follow the same
arguments of step 2 in \cite{CK}. 

Assume
\[
t^2e_{2k-1}\in \frac{1}{(t\lambda)^{2k}}IS^1(R(\log R)^{2k-1},\,\mathcal{Q}'_{k-1})
\]
is given. We begin by isolating the leading component $e_{2k-1}^0$
which includes the terms of top degree in $R$ as well as those of one
degree less (the rest will merge into $e_{2k}$, see step 4 below). Thus we write
$$t^2e^0_{2k-1}=\frac{1}{(t\lambda)^{2k-1}}\sum_{j=0}^{2k-1}aq_j(a)(\log
R)^j+\frac{1}{(t\lambda)^{2k}}\sum_{j=0}^{2k}\tilde{q}_j(a)(\log
R)^j$$
Consider the following equation
$$t^2\widetilde{L}(v_{2k})=t^2e^0_{2k-1}$$
where $\widetilde{L}$ is
$$\widetilde{L}:= -\partial_t^2+\partial^2_r+\frac 1r\partial_r-\frac {1}{r^2}$$

Homogeneity considerations suggest that we should look for a solution
$v_{2k}$ which has the form (notice here we already switched into $R$)
$$v_{2k}=\frac{1}{(t\lambda)^{2k-1}}\sum_{j=0}^{2k-1}W_{2k}^j(a)(\log R)^j+\frac{1}{(t\lambda)^{2k}}\sum_{j=0}^{2k}\widetilde{W}_{2k}^j(a)(\log R)^j$$
The one-dimensional equations for $W^j_{2k}$, $\widetilde{W}^j_{2k}$
are obtained by matching the powers of $\log R$. Then we conjugate out the power
of $t$ and rewrite the systems in the $a$ variable, we get (see step 2 in \cite{CK}
for details)
\begin{align*}
&\mathcal{L}_{(2k-1)\nu}W^j_{2k}=aq_j(a)-F_j(a)\\
&\mathcal{L}_{2k\nu}\widetilde{W}^i_{2k}=\tilde{q}_i(a)-\widetilde{F}_i(a)
\end{align*}
the definition of $\mathcal{L}_{\beta}$ is following \cite{CK}. Solving this system with Cauchy data at $a=0$ yields
solutions which satisfy
\begin{align*}
&W_{2k}^j(a)\in a^3\mathcal{Q}_{k},\quad j=\overline{0,2k-1}\\
&\widetilde{W}_{2k}^i\in a^2\mathcal{Q}_k,\quad i=\overline{0,2k}
\end{align*}
This is guaranteed by lemma 3.9 from \cite{KST0}

To finish this step, we need to make a adjustment for
$v_{2k}$ because of the singularity of $\log R$ at $R=0$. Also, we need to make sure that $v_{2k}$ has order 3 vanishing at $R= 0$. Thus we
define $v_{2k}$ as
\begin{align*}
&v_{2k}:=\\&
\frac{1}{(t\lambda)^{2k-1}}\sum_{j=0}^{2k-1}W^j_{2k}(a)\Big(\frac 12
\log (1+R^2)\Big)^j+\frac{1}{(t\lambda)^{2k}}\frac{R}{(1+R^2)^{\frac{1}{2}}}\sum_{j=0}^{2k}\widetilde{W}^j_{2k}(a)\Big(\frac 12
\log (1+R^2)\Big)^j
\end{align*}
We will get a large error near $R=0$, but it is not very
important since the purpose of the correction is to improve the
error near large $R$. Since $a=R/t\lambda$, it's easy to pull
out a $a^3$ factor from $W$'s and $a^2$ from $\widetilde{W}$'s to see
that we have \eqref{eq:v_2k}.\\

{\bf{Step 4}}: {\it{Show that the error $e_{2k}$ generated by $u_{2k}
    = u_{2k-1}+ v_{2k}$ satisfies \eqref{eq:e_2k}.}} \\

Write 
$$t^2e_{2k}=t^2(e_{2k-1}-e^0_{2k-1})+t^2\big(e^0_{2k-1}-
(-\partial_t^2+\partial_r^2+\frac 1r \partial_r-\frac {1}{r^2})(v_{2k})\big)+t^2N_{2k}(v_{2k})$$
where we recall that except the nonlinear term $t^2N_{2k}(v_{2k})$ the
rest is proved satisfying \eqref{eq:e_2k} following the same arguments
as step 3 in \cite{CK}. For the term $t^2N_{2k}(v_{2k})$, the main method here is to split the nonlinear term in
three parts
\begin{align*}
-t^2N_{2k}(v_{2k})=I+II+III=&a^{-2}\Big[\big(f(u_{2k-1}+v_{2k})-f(u_{2k-2})-f'(u_{2k-1})\big)v_{2k}\Big]\\
+&a^{-2}\Big[\big(f'(u_{2k-1})-f'(u_0)\big)v_{2k}\Big]+a^{-2}\Big[\big(f'(u_0)-1\big)v_{2k}\Big]\\
\end{align*}
and prove each of them lies in a sub-space of what we need in $(\ref{eq:e_2k})$
\begin{align*}
&I\in
a^6\frac{1}{(t\lambda)^{2k}}\sum_{\beta=b,b_{1,2}}\beta IS^1\Big(R(\log
R)^{2k-1},\mathcal Q'_k\Big)\\
&II\in
a^2\frac{1}{(t\lambda)^{2k}}\sum_{\beta=b,b_{1,2}}\beta IS^1\Big(R(\log
R)^{2k-1},\mathcal Q'_k\Big)\\
&III\in
a^2\frac{1}{(t\lambda)^{2k}}IS^3\Big(R^{-1}(\log
R)^{2k},\mathcal Q_k\Big)\\
\end{align*}
The arguments to prove those mimic section 3.8.3 in \cite{catalin}. 
\begin{remark}
One might have doubts since the function space $IS^k(R^m(\log R)^l)$
we are using here is different than \cite{catalin}. To verify this, one just needs to see that the
function spaces defined in \cite{catalin} are the subspaces of our
new defined function space in \cite{CK}. Thus the argument in \cite{catalin} applies to
our case.  
\end{remark}

Iteration of {\bf{Step 1}} - {\bf{Step 4}} immediately furnishes the proof of Theorem~\ref{thm:approxsol} . 

\end{proof}

\section*{Acknowledgement}
The author thanks Joachim Krieger and Willie Wong for many stimulating discussions
and helpful suggestions.

\bigskip


\centerline{\scshape Can Gao }
\medskip
{\footnotesize
 \centerline{B\^{a}timent des Math\'ematiques, EPFL}
\centerline{Station 8,
CH-1015 Lausanne,
  Switzerland}
  \centerline{\email{can.gao@epfl.ch}}
} 

\end{document}